\documentclass[11pt]{article}%
\usepackage{amsmath}
\usepackage{amsfonts}
\usepackage{amssymb}
\usepackage{graphicx}
\usepackage{a4wide}%
\providecommand{\U}[1]{\protect\rule{.1in}{.1in}}
\newtheorem{theorem}{Theorem}

\newtheorem{corollary}[theorem]{Corollary}

\newtheorem{example}[theorem]{Example}

\newtheorem{lemma}[theorem]{Lemma}

\newtheorem{proposition}[theorem]{Proposition}
\newtheorem{remark}[theorem]{Remark}

\begin{document}

\title{On the duality gap and Gale's example in conic linear programming}
\author{Constantin Z\u{a}linescu\thanks{Octav Mayer Institute of Mathematics, Ia\c{s}i
Branch of Romanian Academy, Ia\c{s}i, Romania, and University
\textquotedblleft Alexandru Ioan Cuza" Ia\c{s}i, Romania; email:
\texttt{zalinesc@uaic.ro}.}}
\date{}
\maketitle

\begin{abstract}
The aim of this paper is to revisit some duality results in conic linear
programming and to answer an open problem related to the duality gap function
for Gale's example.

\end{abstract}

\section{Introduction}

In this paper we are mainly concerned by duality results for the following
conic linear programming problem

\medskip(P$_{c,b}$) \ minimize $\ c(x)$ \ s.t. \ $x\in P,\ Ax-b\in Q$,

\medskip\noindent and its dual

\medskip(P$_{c,b}^{\ast}$) \ maximize $\ y^{\ast}(b)$ \ s.t. $\ y^{\ast}\in
Q^{+}$, $A^{\ast}y^{\ast}-c\in-P^{+}$,

\medskip\noindent where $X$, $Y$ are Hausdorff locally convex spaces,
$X^{\ast}$ and $Y^{\ast}$ are their topological dual spaces, $A:X\rightarrow
Y$ is a continuous linear operator, $A^{\ast}:Y^{\ast}\rightarrow X^{\ast}$ is
the adjoint of $A$, $P\subset X$ and $Q\subset Y$ are convex cones,
$P^{+}\subset X^{\ast}$ and $Q^{+}\subset Y^{\ast}$ are the positive dual
cones of $P$ and $Q$, $b\in Y$ and $c\in X^{\ast}$.

\smallskip

In dealing with duality results for the problem (P$_{c,b}$), one uses several
approaches in the literature. For example, in the paper \cite{Kre61}, the book
\cite{AndNas87}, as well as in the recent paper \cite{KhMoTr19}, one
associates certain convex sets to (P$_{c,b}$) and (P$_{c,b}^{\ast}$) and
studies their relationships; in \cite{Sha01} one uses the Lagrangian function
associated to problem (P$_{c,b}$) and several results from convex analysis.

\smallskip In this paper we derive the main duality results for problems
(P$_{c,b}$) and (P$_{c,b} ^{\ast}$) using Rockafellar's perturbation method
(see \cite{Roc74}, \cite{Zal02}). Then, we give an answer to the open problem
concerning the perturbed Gale's example considered in \cite[p.\ 12]%
{ViKiTaYe16}, followed by a discussion of what seems to be the first version
of Gale's example.

\smallskip Below, we introduce some notions, notations and preliminary results.

\smallskip Having $X$ a Hausdorff locally convex space (H.l.c.s.\ for short),
$X^{\ast}$ is its topological dual endowed with its weakly-star topology
$w^{\ast}:=\sigma(X^{\ast},X)$. The value $x^{\ast}(x)$ of $x^{\ast}\in
X^{\ast}$ at $x\in X$ is denoted by $\left\langle x,x^{\ast}\right\rangle $.
It is well known that $(X^{\ast},w^{\ast})^{\ast}$ can be identified with $X$,
what we do in the sequel. For $E\subset X$, by $\operatorname*{span}E$,
$\operatorname*{aff}E$, $\operatorname*{int}E$ and $\operatorname*{cl}E$ one
denotes the \emph{linear hull}, the \emph{affine hull}, the \emph{interior}
and the \emph{closure} of $E$, respectively; moreover, the \emph{intrinsic
core} (or \emph{relative algebraic interior}) of $\emptyset\neq E\subset X$ is
the set
\[
\operatorname*{icr}E:=\{x_{0}\in X\mid\forall x\in\operatorname*{aff}%
E,\ \exists\delta\in\mathbb{P},\ \forall\lambda\in\mathbb{R}:\left\vert
\lambda\right\vert \leq\delta\Rightarrow(1-\lambda)x_{0}+\lambda x\in E\},
\]
where $\mathbb{P}:={}]0,\infty[$; the \emph{core} (or \emph{algebraic
interior}) of $E,$ denoted $\operatorname*{cor}E$, is $\operatorname*{icr}E$
if $\operatorname*{aff}E=X$ and the empty set otherwise. Having $(\emptyset
\neq)$ $K\subset X$ a convex cone (that is, $x+x^{\prime}\in K$ and $tx\in K$
for all $x,x^{\prime}\in K$ and $t\in\mathbb{R}_{+}:=[0,\infty\lbrack$), we
set $x\leq_{K}x^{\prime}$ (equivalently $x^{\prime}\geq_{K}x$) for
$x,x^{\prime}\in X$ with $x^{\prime}-x\in K$; clearly $\leq_{K}$ is a
\emph{preorder} on $X$, that is, $\leq_{K}$ is reflexive and transitive. For
$\emptyset\neq A\subset X$ (and similarly for $\emptyset\neq B\subset X^{\ast
}$) we set $A^{+}:=\{x^{\ast}\in X^{\ast}\mid\forall x\in A:\left\langle
x,x^{\ast}\right\rangle \geq0\}$ for the \emph{positive dual cone} of $A$; it
is well known that $A^{+}$ $(\subset X^{\ast})$ is a $w^{\ast}$-closed convex
cone and $(K^{+})^{+}=\operatorname*{cl}K$ if $K\subset X$ is a convex cone.

\smallskip Having a function $f:X\rightarrow\overline{\mathbb{R}}%
:=\mathbb{R}\cup\{-\infty,\infty\}$, its \emph{domain} is the set
$\operatorname*{dom}f:=\{x\in X\mid f(x)<\infty\}$; $f$ is \emph{proper} if
$\operatorname*{dom}f\neq\emptyset$ and $f(x)\neq-\infty$ for all $x\in X$;
$f$ is \emph{convex} if its \emph{epigraph} $\operatorname*{epi}f:=\{(x,t)\in
X\times\mathbb{R}\mid f(x)\leq t\}$ is convex; $f$ is \emph{positively
homogeneous} if $f(tx)=tf(x)$ for all $t\in\mathbb{P}$ and $x\in X$; $f$ is
\emph{sub\-additive} if $f(x+x^{\prime})\leq f(x)+f(x^{\prime})$ for all
$x,x^{\prime}\in\operatorname*{dom}f$; $f$ is \emph{sublinear} if $f$ is
positively homogeneous, sub\-additive and $f(0)=0$; $f$ is \emph{lower
semi\-continuous} (l.s.c.\ for short) \emph{at} $x\in X$ if $\liminf
_{x^{\prime}\rightarrow x}f(x^{\prime})\geq f(x)$, where $\overline
{\mathbb{R}}$ is endowed with its usual topology and partial order; $f$ is
\emph{l.s.c.}\ if $f$ is l.s.c.\ at any $x\in X$; the \emph{l.s.c.\ envelope}
of $f$ is the function $\overline{f}:X\rightarrow\overline{\mathbb{R}}$ such
that $\operatorname*{epi}\overline{f}=\operatorname*{cl}(\operatorname*{epi}%
f)$, and so $\overline{f}$ is convex if $f$ is so; the
\emph{sub\-differential} of $f$ at $x\in X$ with $f(x)\in\mathbb{R}$ is the
set
\[
\partial f(x):=\{x^{\ast}\in X^{\ast}\mid\forall x^{\prime}\in X:\left\langle
x^{\prime}-x,x^{\ast}\right\rangle \leq f(x^{\prime})-f(x)\}
\]
and $\partial f(x):=\emptyset$ if $f(x)\notin\mathbb{R}$; $f$ is
\emph{sub\-differentiable} at $x\in X$ if $\partial f(x)\neq\emptyset$. The
\emph{conjugate} of $f$ is the function $f^{\ast}:X^{\ast}\rightarrow
\overline{\mathbb{R}}$ defined by
\[
f^{\ast}(x^{\ast}):=\sup\{\left\langle x,x^{\ast}\right\rangle -f(x)\mid x\in
X\}=\sup\{\left\langle x,x^{\ast}\right\rangle -f(x)\mid x\in
\operatorname*{dom}f\}\quad(x^{\ast}\in X^{\ast}),
\]
where $\sup\emptyset:=-\infty$; clearly, $f^{\ast}$ is a $w^{\ast}%
$-l.s.c.\ convex function. Having $g:X^{\ast}\rightarrow\overline{\mathbb{R}}%
$, its conjugate $g^{\ast}:X\rightarrow\overline{\mathbb{R}}$ is defined
similarly. Notice that $f^{\ast}=(\overline{f})^{\ast}$; moreover, for $x\in
X$ and $x^{\ast}\in X^{\ast}$ one has%
\[
x^{\ast}\in\partial f(x)\Leftrightarrow\lbrack f(x)\in\mathbb{R}%
\ \wedge\ f(x)+f^{\ast}(x^{\ast})=\left\langle x,x^{\ast}\right\rangle
]\Rightarrow\overline{f}(x)=f(x)\in\mathbb{R}\Rightarrow\partial
f(x)=\partial\overline{f}(x).\label{r-goz2b}%
\]

As in \cite{Zal02}, the class of proper convex functions defined on $X$ is
denoted by $\Lambda(X)$.

It is worth observing that, for $f:X\rightarrow\overline{\mathbb{R}}$ a convex
function, one has
\begin{equation}
\left[  f(x_{0})=-\infty\ \ \text{and \ }x\in\operatorname*{icr}%
(\operatorname*{dom}f)\right]  \Rightarrow f(x)=-\infty\label{r-goz4}%
\end{equation}
for $x_{0}\in X$ (by \cite[Prop.\ 2.1.4]{Zal02}) and
\begin{gather}
\overline{f}\in\Lambda(X)\Leftrightarrow\left[  \exists x\in X:\overline
{f}(x)\in\mathbb{R}\right]  \Leftrightarrow\left[  \exists x^{\ast}\in
X^{\ast}:f^{\ast}(x^{\ast})\in\mathbb{R}\right]  \Rightarrow f^{\ast\ast
}=\overline{f},\nonumber\\
\left[  \exists x\in X:\overline{f}(x)=-\infty\right]  \Leftrightarrow\left[
\forall x\in\operatorname*{cl}(\operatorname*{dom}f):\overline{f}%
(x)=-\infty\right]  \Leftrightarrow f^{\ast}=\infty\Leftrightarrow f^{\ast
\ast}=-\infty,\nonumber
\end{gather}
by \cite[Ths.\ 2.2.6, 2.3.4]{Zal02}. The \emph{directional derivative} of
$f\in\Lambda(X)$ \emph{at} $x\in\operatorname*{dom}f$ is
\[
f_{+}^{\prime}(x,\cdot):X\rightarrow\overline{\mathbb{R}},\quad f_{+}^{\prime
}(x,u):=\lim_{t\rightarrow0+}\frac{f(x+tu)-f(x)}{t}=\inf_{t>0}\frac
{f(x+tu)-f(x)}{t}.
\]
One has $f_{+}^{\prime}(x,0)=0$, $f_{+}^{\prime}(x,su)=sf_{+}^{\prime}(x,u)$
and $f_{+}^{\prime}(x,u+u^{\prime})\leq f_{+}^{\prime}(x,u)+f_{+}^{\prime
}(x,u^{\prime})$ for $s\in\mathbb{P}$ and $u,u^{\prime}\in X$; hence
$f_{+}^{\prime}(x,\cdot)$ is sublinear. It follows easily that
\[
f_{+}^{\prime}(x,u)\leq f(x+u)-f(x),\quad\partial f(x)=\partial f_{+}^{\prime
}(x,\cdot)(0)\quad\forall f\in\Lambda(X),\ x\in\operatorname*{dom}f,\ u\in
X.\label{r-goz10}%
\]

The \emph{indicator function} of $E\subset X$ is $\iota_{E}:X\rightarrow
\overline{\mathbb{R}}$ defined by $\iota_{E}(x):=0$ for $x\in E$ and
$\iota_{E}(x):=\infty$ for $x\in X\setminus E$; notice that $\iota_{E}$ is
l.s.c.\ iff $E$ is closed, and $\iota_{E}$ is convex iff $E$ is convex.

\section{Alternative proofs for some duality results in conic linear
programming}

In the sequel, we consider the H.l.c.s.\ $X$ and $Y$, the continuous linear
operator $A:X\rightarrow Y$, its adjoint $A^{\ast}:Y^{\ast}\rightarrow
X^{\ast}$ defined by $A^{\ast}y^{\ast}:=y^{\ast}\circ A$ (and so $\left\langle
Ax,y^{\ast}\right\rangle =\left\langle x,A^{\ast}y^{\ast}\right\rangle $ for
$x\in X$, $y^{\ast}\in Y^{\ast}$), the convex cones $P\subset X$, $Q\subset
Y$, as well as their positive dual cones $P^{+}\subset X^{\ast}$ and
$Q^{+}\subset Y^{\ast}$. The preorders defined by $P$, $Q$, $P^{+}$ and
$Q^{+}$ are simply denoted by $\leq$.

For $c\in X^{\ast}$ we associate the mapping $\Phi_{c}:X\times Y\rightarrow
\overline{\mathbb{R}}$, defined by
\begin{equation}
\Phi_{c}(x,y):=\left\langle x,c\right\rangle +\iota_{P}(x)+\iota
_{Q}(Ax-y)>-\infty\quad\big((x,y)\in X\times Y\big); \label{r-goz5}%
\end{equation}
$\Phi_{c}\ $is a proper sublinear function because $\iota_{P}$ and $\iota_{Q}$
are so.

For $(x^{\ast},y^{\ast})\in X^{\ast}\times Y^{\ast}$ one has%
\begin{align*}
\Phi_{c}^{\ast}(x^{\ast},y^{\ast})  &  =\sup\left\{  \left\langle x,x^{\ast
}\right\rangle +\left\langle y,y^{\ast}\right\rangle -\left\langle x,c^{\ast
}\right\rangle -\iota_{P}(x)-\iota_{Q}(Ax-y)\mid x\in X,\ y\in Y\right\} \\
&  =\sup\left\{  \left\langle x,x^{\ast}-c^{\ast}\right\rangle +\left\langle
Ax-q,y^{\ast}\right\rangle \mid x\in P,\ q\in Q\right\} \\
&  =\sup\left\{  \left\langle x,x^{\ast}-c^{\ast}+A^{\ast}y^{\ast
}\right\rangle -\left\langle q,y^{\ast}\right\rangle \mid x\in P,\ q\in
Q\right\} \\
&  =\left\{
\begin{array}
[c]{ll}%
0 & \text{if }y^{\ast}\in Q^{+}\text{ and }c^{\ast}-A^{\ast}y^{\ast}-x^{\ast
}\in P^{+},\\
\infty & \text{otherwise,}%
\end{array}
\right. \\
&  =\iota_{P^{+}}(c^{\ast}-A^{\ast}y^{\ast}-x^{\ast})+\iota_{Q^{+}}(y^{\ast}).
\end{align*}

The \emph{marginal} (\emph{value}) \emph{function} associated to $\Phi_{c}$
is
\begin{equation}
h_{c}:Y\rightarrow\overline{\mathbb{R}},\quad h_{c}(y):=\inf\{\Phi
_{c}(x,y)\mid x\in X\}=\inf\{\left\langle x,c\right\rangle \mid x\geq
0,\ Ax\geq y\},\nonumber
\end{equation}
where $\inf\emptyset:=\infty$; hence $\operatorname*{dom}h_{c}=A(P)-Q$.

Clearly, $h_{c}(0)\leq\Phi_{c}(0,0)=0$ and $h_{c}$ is increasing, that is,
$\ h_{c}(y_{1})\leq h_{c}(y_{2})$ for $y_{1}\leq y_{2}$. Because $\Phi_{c}$ is
convex, so is $h_{c}$. Moreover, for $\alpha>0$ and $y\in Y$ one has
\[
h_{c}(\alpha y)=\inf\nolimits_{x\in X}\Phi_{c}(x,\alpha y)=\inf\nolimits_{x\in
X}\alpha\Phi_{c}(\alpha^{-1}x,y)=\alpha\inf\nolimits_{x^{\prime}\in X}\Phi
_{c}(x^{\prime},y)=\alpha h_{c}(y)\text{,}%
\]
and so $h_{c}$ is positively homogeneous; it follows that $h_{c}%
(0)\in\{0,-\infty\}$.

If $h_{c}(0)=-\infty$, then $h_{c}(y)=h_{c}(y+0)\leq h_{c}(y)+h_{c}%
(0)=-\infty$ for every $y\in\operatorname*{dom}h_{c}$; moreover, $\partial
h_{c}(0)=\emptyset$ (by definition), $h_{c}^{\ast}=\infty$ and $h_{c}%
^{\ast\ast}=-\infty$.

Assume that $h_{c}(0)=0;$ then%
\begin{align*}
y^{\ast}\in\partial h_{c}(0)  &  \Leftrightarrow\lbrack\forall y\in
\operatorname*{dom}h_{c}:\left\langle y,y^{\ast}\right\rangle \leq
h_{c}(y)]\Leftrightarrow\lbrack\forall x\in P,\ \forall z\in Q:\left\langle
Ax-z,y^{\ast}\right\rangle \leq\left\langle x,c\right\rangle ]\\
&  \Leftrightarrow\left[  y^{\ast}\in Q^{+}\wedge\lbrack\forall x\in
P:0\leq\left\langle x,c-A^{\ast}y^{\ast}\right\rangle ]\right]
\Leftrightarrow y^{\ast}\in Q^{+}\cap(A^{\ast})^{-1}(c-P^{+}),
\end{align*}
and so $\partial h_{c}(0)=Q^{+}\cap(A^{\ast})^{-1}(c-P^{+})$.

It follows easily from the definition of $h_{c}$ (or from \cite[Th.\ 2.6.1(i)]%
{Zal02}) that%
\begin{align}
h_{c}^{\ast}(y^{\ast})  &  =\Phi_{c}^{\ast}(0,y^{\ast})=\left\{
\begin{array}
[c]{ll}%
0 & \text{if }y^{\ast}\in Q^{+}\text{ and }A^{\ast}y^{\ast}-c\in(-P^{+}),\\
\infty & \text{otherwise,}%
\end{array}
\right. \nonumber\\
&  =\iota_{P^{+}}(c-A^{\ast}y^{\ast})+i_{Q^{+}}(y^{\ast})=i_{Q^{+}}(y^{\ast
})+\iota_{-P^{+}}(A^{\ast}y^{\ast}-c), \label{r-g13}%
\end{align}
and so $\operatorname*{dom}h_{c}^{\ast}=Q^{+}\cap(A^{\ast})^{-1}(c-P^{+})$.
Consequently, if $h_{c}^{\ast}$ is proper (hence $h_{c}(0)=0$), then%
\begin{equation}
\partial h_{c}(0)=\operatorname*{dom}h_{c}^{\ast}=Q^{+}\cap(A^{\ast}%
)^{-1}(c-P^{+})\neq\emptyset,\ \ h_{c}^{\ast}=\iota_{\partial h_{c}(0)}\text{
\ and \ }\overline{h_{c}}=h_{c}^{\ast\ast}. \label{r-goz3}%
\end{equation}

On the other hand,
\[
h_{c}(0)=0\Leftrightarrow\left[  x\in P\cap A^{-1}(Q)\Rightarrow\left\langle
x,c\right\rangle \geq0\right]  \Leftrightarrow c\in\lbrack P\cap
A^{-1}(Q)]^{+}.
\]
If $P$ and $Q$ are closed, then $\left(  A^{-1}(Q)\right)  ^{+}%
=\operatorname*{cl}_{w^{\ast}}A^{\ast}(Q^{+})$ by \cite[Lem.\ 4]{Kre61} (see
also \cite[Lem.\ 1]{Zal78}), and so%
\[
\lbrack P\cap A^{-1}(Q)]^{+}=\operatorname*{cl}\nolimits_{w^{\ast}}%
\big[P^{+}+\left(  A^{-1}(Q)\right)  ^{+}\big]=\operatorname*{cl}%
\nolimits_{w^{\ast}}\left[  P^{+}+A^{\ast}(Q^{+})\right]  .
\]

For $c\in X^{\ast}$ and $b\in Y$, consider
\begin{equation}
\Phi_{c,b}:X\times Y\rightarrow\overline{\mathbb{R}},\quad\Phi_{c,b}%
(x,y):=\Phi_{c}(x,y+b); \label{r-g11a}%
\end{equation}
hence%
\begin{equation}
\forall y\in Y:h_{c,b}(y)=h_{c}(y+b)\ \wedge\ \overline{h_{c,b}}%
(y)=\overline{h_{c}}(y+b)\ \wedge\ \partial h_{c,b}(y)=\partial h_{c}(y+b),
\label{r-g11c}%
\end{equation}
where $h_{c,b}$ is the value function associated to $\Phi_{c,b}$. It follows
that $\Phi_{c,b}^{\ast}(x^{\ast},y^{\ast})=\Phi_{c}^{\ast}(x^{\ast},y^{\ast
})-\left\langle b,y^{\ast}\right\rangle $ for $(x^{\ast},y^{\ast})\in X^{\ast
}\times Y^{\ast}$, whence
\begin{equation}
\forall y^{\ast}\in Y^{\ast}:h_{c,b}^{\ast}(y^{\ast})=\Phi_{c,b}^{\ast
}(0,y^{\ast})=\Phi_{c}^{\ast}(0,y^{\ast})-\left\langle b,y^{\ast}\right\rangle
=h_{c}^{\ast}(y^{\ast})-\left\langle b,y^{\ast}\right\rangle . \label{r-g11b}%
\end{equation}

From (\ref{r-goz5}) and (\ref{r-g11a}), one has $\Phi_{c,b}(x,0)=\Phi
_{c}(x,b)=\left\langle x,c\right\rangle $ if $Ax\geq b$, $x\geq0$, and
$\Phi_{c,b}(x,0)=\infty$ otherwise, and so the problem (P$_{c,b}$) becomes

\medskip(P$_{c,b}$) \ minimize \ $\Phi_{c,b}(x,0)$ \ s.t.\ $x\in X$.

\medskip On the other hand, from (\ref{r-g13}) and (\ref{r-g11b}) one has
$\Phi_{c,b}^{\ast}(0,y^{\ast})=-\left\langle b,y^{\ast}\right\rangle $ if
$A^{\ast}y^{\ast}\leq c$, $y^{\ast}\geq0$, and $\Phi_{c,b}^{\ast}(0,y^{\ast
})=\infty$ otherwise, and so the problem (P$_{c,b}^{\ast}$) becomes

\medskip(P$_{c,b}^{\ast}$) \ maximize \ $-\Phi_{c,b}^{\ast}(0,y^{\ast})$
\ s.t.\ $y^{\ast}\in Y^{\ast}$.

\medskip For further use, similarly to the notations from \cite{ViKiTaYe16},
we set:

\smallskip$F(b):=\{x\in X\mid Ax\geq b,\ x\geq0\}$ -- the feasible set of
(P$_{c,b}$);

\smallskip$\varphi(c,b):=\inf\{\left\langle x,c\right\rangle \mid x\in
F(b)\}=h_{c}(b)\in\overline{\mathbb{R}}$ -- the value of (P$_{c,b}$);

\smallskip$F^{\ast}(c):=\{y^{\ast}\in Y^{\ast}\mid A^{\ast}y^{\ast}\leq
c,\ y^{\ast}\geq0\}$ -- the feasible set of (P$_{c,b}^{\ast}$);

\smallskip$\psi(c,b):=\sup\{\left\langle b,y^{\ast}\right\rangle \mid y^{\ast
}\in F^{\ast}(c)\}=h_{c,b}^{\ast\ast}(b)\in\overline{\mathbb{R}}$ -- the value
of (P$_{c,b}^{\ast}$);

\smallskip$\Lambda:=\{(c,b)\in X^{\ast}\times Y\mid F^{\ast}(c)\neq
\emptyset,\ F(b)\neq\emptyset\}=(P^{+}+A^{\ast}(Q^{+}))\times(A(P)-Q)$;

\smallskip$g:\Lambda\rightarrow\mathbb{R}$, $g(c,b):=\varphi(c,b)-\psi(c,b)$
-- the duality gap function of (P$_{c,b}$).

\medskip

Applying \cite[Th.\ 2.6.1]{Zal02} for $\Phi_{c,b}$ one gets the following result.

\begin{proposition}
\label{p4}The following assertions hold:

\emph{(i)} $h_{c}(b)=$\emph{val(P}$_{c,b}$\emph{)}; $h_{c}^{\ast\ast}%
(b)=$\emph{val(P}$_{c,b}^{\ast}$\emph{)}; \emph{val(P}$_{c,b}$\emph{)}$\geq
$\emph{val(P}$_{c,b}^{\ast}$\emph{)} for all $c\in X^{\ast}$ and $b\in Y$;
\emph{val(P}$_{c,b}$\emph{)}$<\infty$ $\Leftrightarrow$ $b\in A(P)-Q$;
\emph{val(P}$_{c,b}^{\ast}$\emph{)}$>-\infty$ $\Leftrightarrow$ $c\in
P^{+}+A^{\ast}(Q^{+})$; consequently, $g(c,b)\in\mathbb{R}_{+}$ for all
$(c,b)\in\Lambda$.

\emph{(ii)} [$h_{c}(b)\in\mathbb{R}$ and $h_{c}$ is l.s.c.\ at $b$]
$\Leftrightarrow$ \emph{val(P}$_{c,b}$\emph{)}$=$\emph{val(P}$_{c,b}^{\ast}%
$\emph{)}$\in\mathbb{R}$.

\emph{(iii)} $\partial h_{c}(b)\neq\emptyset\Leftrightarrow$ [\emph{val(P}%
$_{c,b}$\emph{)}$=$\emph{val(P}$_{c,b}^{\ast}$\emph{)}$\in\mathbb{R}$ and
\emph{(P}$_{c,b}^{\ast}$\emph{)} has optimal solutions]; moreover,
\emph{Sol(P}$_{c,b}^{\ast}$\emph{)}$=\partial h_{c}(b)$ if $\partial
h_{c}(b)\neq\emptyset$, where \emph{Sol(P}$_{c,b}^{\ast}$\emph{)} is the set
of optimal solutions of problem \emph{(P}$_{c,b}^{\ast}$\emph{)}.

\emph{(iv)} $\overline{h_{c}}$ is proper~$\Leftrightarrow h_{c}^{\ast}$ is
proper~$\Leftrightarrow$ $\operatorname*{dom}h_{c}^{\ast}\neq\emptyset
~\Leftrightarrow$ $c\in P^{+}+A^{\ast}(Q^{+})$.
\end{proposition}

Proof. We apply \cite[Th.\ 2.6.1]{Zal02} for $\Phi_{c,b}$ and its associated
value function $h_{c,b}$. Having in view the equality $h_{c,b}(\cdot
)=h_{c}(\cdot+b)$ from (\ref{r-g11c}), one has $h_{c,b}^{\ast\ast}%
(\cdot)=h_{c}^{\ast\ast}(\cdot+b)$, too. Hence
\begin{equation}
h_{c,b}(0)=h_{c}(b),\quad h_{c,b}^{\ast\ast}(0)=h_{c}^{\ast\ast}%
(b),\quad\partial h_{c,b}(0)=\partial h_{c}(b),\quad h_{c,b}^{\ast}%
=h_{c}^{\ast}-\left\langle b,\cdot\right\rangle . \label{r-g14}%
\end{equation}

(i) Using (\ref{r-g14}) and \cite[Th.\ 2.6.1(iii)]{Zal02} one gets
\[
\text{val(}P_{c,b}\text{)}=h_{c,b}(0)=h_{c}(b)\geq h_{c}^{\ast\ast}%
(b)=h_{c,b}^{\ast\ast}(0)=\text{val(}P_{c,b}\text{)}.
\]
The conclusion follows observing that [val(P$_{c,b}$)$<\infty$
$\Leftrightarrow$ $F(b)\neq\emptyset$] and $b$[val(P$_{c,b}^{\ast}$)$>-\infty
$~$\Leftrightarrow$ $F^{\ast}(c)\neq\emptyset$].

(ii) and (iii) Having in view (\ref{r-g14}), these assertions are nothing else
than assertions (v) and (vi) from \cite[Th.\ 2.6.1]{Zal02} in the present
case, respectively.

(iv) Having in view (\ref{r-g14}), the first equivalence is provided by
\cite[Th.\ 2.6.1(vii)]{Zal02}, while the other two are immediate from
(\ref{r-goz3}). \hfill$\square$

\medskip

Notice that the equality $h_{c}^{\ast\ast}(b)=$val(P$_{c,b}^{\ast}$) from
Proposition \ref{p4}(i) is established in \cite[Prop.\ 2.2]{Sha01}, the other
assertions from (i) being well known; moreover, the assertion (iii) is
established in \cite[Prop.\ 2.5]{Sha01}, while the equivalence in (iii) is
also established in \cite[Th.\ 1]{GrOsZa02}.

\medskip The next result provides sufficient interiority conditions for the
sub\-differentiability of $h_{c}$ at $b$.

\begin{proposition}
\label{p3}Let $Y_{0}:=\operatorname*{span}(A(P)-Q)$ be endowed with the
induced topology. Assume that one of the following conditions is verified:

\emph{(i)} there exists $\lambda_{0}\in\mathbb{R}$ such that $N_{0}:=\{y\in
Y\mid\exists x\in P:Ax-y\in Q$ and $\left\langle x,c\right\rangle \leq
\lambda_{0}\}$ is a neighborhood of $b$ in $Y_{0}$;

\emph{(ii)} there exists $x_{0}\in P$ such that $Ax_{0}-b\in
\operatorname*{int}_{Y_{0}}Q$;

\emph{(iii)} $X$ and $Y$ are Fr\'{e}chet spaces, $P$, $Q$ and $Y_{0}$ are
closed, and $b\in\operatorname*{icr}(A(P)-Q)$;

\emph{(iv)} $\dim Y_{0}<\infty$ and $b\in\operatorname*{icr}(A(P)-Q)$.

Then either \emph{(a)} $h_{c}(y)=-\infty$ for every $y\in\operatorname*{icr}%
(A(P)-Q)$, and so \emph{val(P}$_{c,b}$\emph{)}$=$\emph{val(P}$_{c,b}^{\ast}%
$\emph{)}$=-\infty$, or \emph{(b)}~ $c\in P^{+}+A^{\ast}(Q^{+})$ and
$h_{c}|_{Y_{0}}$ is continuous at $b$. Consequently, in the case \emph{(b)},
$\partial h_{c}(b)\neq\emptyset$, whence \emph{val(P}$_{c,b}$\emph{)}%
$=$\emph{val(P}$_{c,b}^{\ast}$\emph{)}$\in\mathbb{R}$ and \emph{(P}%
$_{c,b}^{\ast}$\emph{)} has optimal solutions; moreover, if $Y_{0}=Y$, then
\emph{Sol(P}$_{c,b}^{\ast}$\emph{)} is $w^{\ast}$-compact.
\end{proposition}

Proof. We apply \cite[Th.\ 2.7.1]{Zal02} for $\Phi:=\Phi_{c,b}$ and its
associated value function $h_{c,b}$.

Clearly, $\operatorname*{dom}\Phi_{c,b}=\cup_{x\in P}\left(  \{x\}\times
(Ax-Q-b)\right)  $, and so $\operatorname*{dom}h_{c,b}=\Pr_{Y}\left(
\operatorname*{dom}\Phi_{c,b}\right)  =A(P)-Q-b=\operatorname*{dom}h_{c}-b$.
Observe that each of the conditions (i)--(iv) implies that $b\in
\operatorname*{icr}(A(P)-Q)=\operatorname*{dom}h_{c}-b$, and so $h_{c}%
(b)<\infty$.

Because $b\in A(P)-Q$ and $P$, $Q$ are convex cones, it follows that
\[
\operatorname*{span}\big(  \Pr\nolimits_{Y}(\operatorname*{dom}\Phi
_{c,b})\big)  =\operatorname*{span}\left(  A(P)-Q\right)  =Y_{0}%
=A(P-P)+Q-Q\supseteq Q\cup A(P)\cup\{b\}.
\]

Observe that, if (i), (ii), (iii) or (iv) is verified, then condition (i),
(iii), (vii) or (viii) from \cite[Th.\ 2.7.1]{Zal02} is verified,
respectively; using this theorem, one obtains that either $h_{c,b}(0)=-\infty
$, or $h_{c,b}(0)\in\mathbb{R}$ and $h_{c,b}|_{Y_{0}}$ is continuous at $0$,
or, equivalently, either $h_{c}(b)=-\infty$, or $h_{c}(b)\in\mathbb{R}$ and
$h_{c}|_{Y_{0}}$ is continuous at $b$.

If $h_{c}(b)=-\infty$, then $h_{c}(y)=-\infty$ for every $y\in
\operatorname*{icr}(\operatorname*{dom}h_{c})$ $[=\operatorname*{icr}%
(A(P)-Q)]$ and val(P$_{c,b}$)$=$val(P$_{c,b}^{\ast}$)$=-\infty$ by Proposition
\ref{p4}(i); hence (a) holds.

If $h_{c}(b)\in\mathbb{R}$, then $h_{c}|_{Y_{0}}$ is proper, whence $\partial
h_{c}(b)\neq\emptyset$ by \cite[Th.\ 2.4.12]{Zal02} (that is (b) holds), and
so val(P$_{c,b}$)$=$val(P$_{c,b}^{\ast}$)\,$\in\mathbb{R}$ and (P$_{c,b}%
^{\ast}$) has optimal solutions by Proposition \ref{p4}(iii). Moreover, if
$Y_{0}=Y$ then $h_{c}$ is continuous at $b$, and so $\partial h_{c}(b)$ is
$w^{\ast}$-compact by \cite[Th.\ 2.4.9]{Zal02}. \hfill$\square$

\smallskip

Observe that Kretschmer obtained, in \cite[Th.\ 3]{Kre61}, that $g(c,b)=0$ and
that (P$_{c,b}^{\ast}$) has optimal solutions when $Y_{0}=Y$, $P,Q$ are
closed, val(P$_{c,b}$)$\in\mathbb{R}$, and condition (iii) above holds; the
version of \cite[Th.\ 3]{Kre61} for the dual problem is stated in
\cite[Cor.\ 3.1]{Kre61}. Notice also that the equality $g(c,b)=0$ is
established for $Y_{0}=Y$ and val(P$_{c,b}$)$\in\mathbb{R}$ in Theorems
3.11--3.13 from \cite{AndNas87}; more precisely, one asks: (iii) above in
\cite[Th.\ 3.13]{AndNas87}; $Q=\{0\}$ and (i) holds in \cite[Th.\ 3.11]%
{AndNas87}; $Q=\{0\}$, $X,Y$ are Banach spaces and $\operatorname*{int}%
P\neq\emptyset$ in \cite[Th.\ 3.12]{AndNas87}. In fact, in the latter case it
is shown that (i) holds. Moreover, the version of \cite[Th.\ 3.13]{AndNas87}
for the dual problem, when $P$ and $Q$ are closed, is stated in
\cite[Cor.\ 3.14]{AndNas87}. In \cite[Sect.\ 3.8]{AndNas87} one finds more
references concerning the absence of duality gap in linear programming.
Furthermore, \cite[Prop.\ 2.9]{Sha01} states that val(P$_{c,b}$)$=$%
val(P$_{c,b}^{\ast}$) and Sol(P$_{c,b}^{\ast}$) is nonempty and bounded
provided that $X,Y$ are Banach spaces and $b\in\operatorname*{int}\left(
A(P)-Q\right)  $; of course, $Y_{0}=Y$ in this case.

\medskip In the last 20 years, the interest for studying optimization problems
in the algebraic framework increased. In fact, several results established in
this framework can be deduced from the corresponding topological ones; the
next corollary is a sample. Actually, having a real linear space $E$, the
strongest locally convex topology on $E$, denoted $\tau_{c}$ and called the
\emph{convex core topology}, is generated by the family of all the semi-norms
defined on $E$ (see, e.g.,\ \cite[Exer.\ 2.10]{Hol75}); notice that the family
$\{N\subset E\mid N$ is convex and $\operatorname*{cor}N\neq\emptyset\}$ is a
neighborhood base of $0\in E$ for $\tau_{c}$. If $\dim E<\infty$, $\tau_{c}$
coincides with any Hausdorff linear topology on $E$. Moreover, $(E,\tau
_{c})^{\ast}=E^{\prime}$, where $E^{\prime}$ is the algebraic dual of $E$, and
so $(E^{\prime},\sigma(E^{\prime},E))^{*}=E$; hence $\tau_{c}$ coincides with
the Mackey topology of $E$ for the dual pair $(E,E^{\prime})$. Because
$(E^{\prime},\tau_{c})^{\ast}=(E^{\prime})^{\prime}$, one has $w^{\ast}%
\neq\tau_{c}$ if $\dim E=\infty$. Having another real linear space $F$, any
linear operator $B:E\rightarrow F$ is continuous when $E$ and $F$ are endowed
with their convex core topologies. Proposition 6.3.1 from \cite{KhaTamZal15b}
collects several results concerning $\tau_{c}$.

\begin{corollary}
\label{c2}Let the real linear spaces $X$ and $Y$ be endowed with their convex
core topologies. Assume that $b\in\operatorname*{icr}(A(P)-Q)$. Then the
conclusions of Proposition \ref{p3} hold.
\end{corollary}

Proof. In the proof of Proposition \ref{p3} we observed that
$\operatorname*{dom}h_{c}=A(P)-Q$, $h_{c}(y)=-\infty$ for every $y\in
\operatorname*{icr}(\operatorname*{dom}h_{c})$ if $h_{c}(b)=-\infty$, and that
$h_{c}|_{Y_{0}}$ is proper if $h_{c}(b)\in\mathbb{R}$, whence $h_{c}$ is
sub\-differentiable on $\operatorname*{icr}(\operatorname*{dom}h_{c})$ by
\cite[Prop.\ 6.3.1(v)]{KhaTamZal15b}. \hfill$\square$

\medskip

In the case $Q:=\{0\}$, Corollary \ref{c2} is comparable with \cite[Th.\ 4]%
{Pin11}. More precisely, the main part of \cite[Th.\ 4]{Pin11}, stated in the
proof of Point 1., establishes (in our terms) that val(P$_{c,b}$%
)$=$val(P$_{c,b}^{\ast}$) whenever $\operatorname*{cor}P\neq\emptyset$ (with
$P:=C$), $Q:=\{0\}$ and val(P$_{c,b}$)$\in\mathbb{R}$. Observe that Corollary
\ref{c2} may not be used for $b\in A(P)\setminus\operatorname*{icr}A(P)$ even
if $\operatorname*{cor}P\neq\emptyset$. Based on the following example, we
asserts that the conclusion of \cite[Th.\ 4]{Pin11} might be false for $b\in
A(P)\setminus\operatorname*{icr}A(P)$. The next example can be found in
\cite{Zal22}.

\begin{example}
\label{ex-z} \emph{(See \cite[Examp.\ 4]{Zal22}.)} Consider $X_{0}%
:=\mathbb{R}^{2}$, $Y:=\mathbb{R}^{3}$, $A_{0}:X_{0}\rightarrow Y$ with
$A_{0}(x_{1},x_{2}):=(x_{1},x_{2},0)$,
\[
P_{0}:=\mathbb{R}\times\mathbb{R}_{+},\quad Q_{0}:=\left\{  (y_{1},y_{2}%
,y_{3})\in Y\mid y_{1},y_{3}\in\mathbb{R}_{+},\ (y_{2})^{2}\leq2y_{1}%
y_{3}\right\}  ,
\]
the conic linear programming problem

\smallskip\emph{(P}$_{y}$\emph{) \ }minimize\ $x_{2}$ \ s.t. $x:=(x_{1}%
,x_{2})\in P_{0}$, $A_{0}x-y\in Q_{0}$,

\smallskip\noindent and $h_{0}:Y\rightarrow\overline{\mathbb{R}}$ defined by
$h_{0}(y):=$ \emph{val(P}$_{y}$\emph{)}. Then $h_{0}(y)=y_{2}$ if
$y:=(y_{1},y_{2},y_{3})\in(\mathbb{R}\times\mathbb{R}_{+}\times\{0\})$,
$h_{0}(y)=0$ for $y\in\mathbb{R}\times\mathbb{R}\times(-\mathbb{P})$, and
$h_{0}(y)=\infty$ elsewhere.
\end{example}

Observe that the problem (P$_{y}$) is equivalent with the following one:

\smallskip(P$_{y}^{\prime}$) \ minimize\ $x_{2}$ \ s.t. $(x,v)\in
P:=P_{0}\times Q_{0}$, $A(x,v)=y$,

\smallskip\noindent where $A:X:=X_{0}\times Y\rightarrow Y$ is defined by
$A(x,v):=A_{0}x-v$. Clearly, the value function associated to problem
(P$_{y}^{\prime}$) is $h_{0}$. By Proposition \ref{p4}(i) one has that
val(P$_{c,b}$)$=h_{0}(b)=b_{2}>0=h_{0}^{\ast\ast}(b)=$val(P$_{c,b}^{\ast}$)
for $b\in\mathbb{R}\times\mathbb{R}\times(-\mathbb{P})=A(P)\setminus
\operatorname*{icr}A(P)$. \hfill$\square$

\begin{corollary}
\label{c1}Let $(c,b)\in\Lambda$ and set $Y_{0}:=\operatorname*{span}(A(P)-Q)$.
Assume that one of the conditions \emph{(i)--(iv)} of Proposition \ref{p3}
holds. Then $\operatorname*{int}_{Y_{0}}(\operatorname*{dom}h_{c}%
)=\operatorname*{int}_{Y_{0}}(A(P)-Q)\neq\emptyset$ and $h_{c}|_{Y_{0}}$ is
finite and continuous at every $y\in\operatorname*{int}_{Y_{0}}%
(\operatorname*{dom}h_{c})$; consequently, all the other conclusions of
Proposition \ref{p3} hold at $y$. Moreover, assume that $y\in
\operatorname*{dom}h_{c}\setminus\operatorname*{int}_{Y_{0}}%
(\operatorname*{dom}h_{c})$ is such that $(h_{c})_{+}^{\prime}(y,\cdot)$ is
proper; then $\partial h_{c}(y)\neq\emptyset$, and so \emph{val(P}$_{c,y}%
$\emph{)}$=$\emph{val(P}$_{c,y}^{\ast}$\emph{)} and \emph{Sol(P}$_{c,y}^{\ast
}$\emph{)}$=\partial h_{c}(y)$.
\end{corollary}

Proof. Recall that having a proper convex function $f:Y\rightarrow
\overline{\mathbb{R}}$ which is continuous at some $y_{0}\in
\operatorname*{dom}f$ (or, equivalently, $f$ is bounded above on some nonempty
open subset of $\operatorname*{dom}f$), then $\operatorname*{int}%
(\operatorname*{dom}f)\neq\emptyset$ and $f$ is continuous at any
$y\in\operatorname*{int}(\operatorname*{dom}f)$ (see, e.g., \cite[Th.\ 2.2.9]%
{Zal02}).

Take $b^{\prime}\in\operatorname*{int}_{Y_{0}}(\operatorname*{dom}h_{c})$;
because, in the present situation, $h_{c}|_{Y_{0}}$ is continuous at $b$, it
follows that $h_{c}|_{Y_{0}}$ is continuous at $b^{\prime}$, and so the other
conclusions of Proposition \ref{p3} hold at $b^{\prime}$.

Assume now that $\overline{b}\in\operatorname*{dom}h_{c}\setminus
\operatorname*{int}_{Y_{0}}(\operatorname*{dom}h_{c})$ is such that
$(h_{c})_{+}^{\prime}(\overline{b},\cdot)$ is proper. Consider $f:=h_{c}%
|_{Y_{0}}$; of course, $f$ is finite and continuous at $b$ and $f_{+}^{\prime
}(\overline{b},\cdot)$ $\left(  =(h_{c})_{+}^{\prime}(\overline{b}%
,\cdot)|_{Y_{0}}\right)  $ is proper. Applying \cite[Prop.\ 2.134(ii)]%
{BonSha00} we get $\left(  \partial h_{c}|_{Y_{0}}(\overline{b})=\right)  $
$\partial f(\overline{b})\neq\emptyset$, whence $\partial h_{c}(\overline
{b})\neq\emptyset$. \hfill$\square$

\begin{remark}
\label{rem1}Notice that, besides the fact that $g(c,b)=0$ for $(c,b)\in\left(
P^{+}+A^{\ast}(Q^{+})\right)  \times\operatorname*{icr}\left(  A(P)-Q\right)
$, Corollary \ref{c1} shows that \emph{Sol(P}$_{c,b}^{\ast}$\emph{)} is
nonempty for such $(c,b)$.
\end{remark}

Observe that Shapiro \cite[Prop.\ 2.8]{Sha01}, when $h_{c}(b)\in\mathbb{R}$
and either $\dim Y<\infty$ or $h_{c}$ is bounded above on a nonempty open set,
established that [val(P$_{c,b}$)$=$val(P$_{c,b}^{\ast}$) and Sol(P$_{c,b}%
^{\ast}$)$\neq\emptyset$] iff $(h_{c})_{+}^{\prime}(b,\cdot)$ is proper (that
is, $\partial h_{c}(b)\neq\emptyset\Leftrightarrow$ $(h_{c})_{+}^{\prime
}(b,\cdot)$ is proper).

Notice also that one uses the same hypotheses on $A$ and $b$ as in Theorems
3.11 and 3.12 from \cite{AndNas87} for getting $g(c,b)=0$ for $(c,b)\in\left(
P^{+}+A^{\ast}(Q^{+})\right)  \times\operatorname*{cor}\left(  A(P)-Q\right)
$ in Theorems 3.1 and 3.2 from \cite{ViKiTaYe16}, respectively; having in view
the discussion after the proof of Proposition \ref{p3}, the conclusions of
\cite[Ths.\ 3.1, 3.2]{ViKiTaYe16} follow by using Corollary \ref{c1}.

\section{On Gale's example in conic linear programming}

We have seen (see Example \ref{ex-z}) that, even in finite-dimensional conic
linear programming, it is possible to have problems with positive duality gap;
of course such examples exist also in the infinite-dimensional case. The main
aim of this section is to answer the open problem from \cite[p.\ 12]%
{ViKiTaYe16} concerning the perturbed Gale's example.

\smallskip Gale's example is presented in Section 3.4.1 of the book
\cite{AndNas87} as being the problem:

\smallskip(PG) minimize $x_{0}$ \ s.t.\ \ $x_{0}+\sum_{k=1}^{\infty}kx_{k}=1$,
$\sum_{k=1}^{\infty}x_{k}=0$, $x_{k}\geq0$, $k=0,1,2$, \ldots .

\medskip In this example, the primal variable space is
\[
X:=\mathbb{R}^{(\mathbb{N})}:=\{(x_{k})_{k\in\mathbb{N}}\subset\mathbb{R}%
\mid\exists k_{0}\in\mathbb{N},\ \forall k\geq k_{0}:x_{k}=0\},
\]
called the generalized finite sequence space; in the sequel $\mathbb{N}%
:=\{0,1,2,...\}$. The dual of (PG) is

\smallskip(PG$^{\ast}$)\ maximize $y_{1}$ \ s.t.\ \ $y_{1}\leq1$,
$ky_{1}+y_{2}\leq0$, $k=1,2,...$,

\smallskip\noindent in which the dual variable space is $\mathbb{R}^{2}$.

As mentioned in \cite{AndNas87}, the only feasible (and so, optimal) solution
of (PG) is $x=(1,0,0,...)$; hence val(PG)$=1$. Moreover, $(y_{1},y_{2}%
)\in\mathbb{R}^{2}$ is feasible for (PG$^{\ast}$) iff $y_{1}\leq0$,
$y_{1}+y_{2}\leq0$, and so $(0,0)$ is an optimal solution of (PG$^{\ast}$).
Hence val(PG$^{\ast}$)$=0$, and so there is a positive duality gap.

\medskip

Observing that $(e_{n})_{n\in\mathbb{N}}$ is an algebraic (Hamel) basis of
$X$, where $e_{n}:=(\delta_{nk})_{k\in\mathbb{N}}$ ($\delta_{nk}$ being the
Kronecker symbol), the algebraic dual $X^{\prime}$ of $X$ is $\mathbb{R}%
^{\mathbb{N}}$, the space of all real sequences $x^{\prime}:=(x_{k}^{\prime
})_{k\geq0}$, with
\[
x^{\prime}(x):=\left\langle x,x^{\prime}\right\rangle :=\sum\nolimits_{k=0}%
^{\infty}x_{k}x_{k}^{\prime}\quad\left(  x:=(x_{k})_{k\in\mathbb{N}}\in
X\right)  .\label{r-g1}%
\]

Take $c:=e_{0}\in X^{\prime}$, $b:=(1,0)\in\mathbb{R}^{2}$,
\begin{gather}
P:=\mathbb{R}_{+}^{(\mathbb{N})}:=\{x\in\mathbb{R}^{(\mathbb{N})}\mid\forall
k\in\mathbb{N}:x_{k}\geq0\},\ \ Q:=\{0\}\subset\mathbb{R}^{2},\label{r-g2a}\\
A:\mathbb{R}^{(\mathbb{N})}\rightarrow\mathbb{R}^{2},\quad A(x):=\left(
x_{0}+\sum\nolimits_{k=1}^{\infty}kx_{k},\sum\nolimits_{k=1}^{\infty}%
x_{k}\right)  ; \label{r-g2b}%
\end{gather}
then $P$ and $Q$ are convex cones and $A$ is a linear operator

We endow $Y:=\mathbb{R}^{2}$ with its usual inner product, and so $Y^{\ast}$
is identified with $Y$; moreover, we endow $X:=\mathbb{R}^{(\mathbb{N})}$ with
the convex core topology $\tau_{c}$, and so $X^{\prime}=\mathbb{R}%
^{\mathbb{N}}$ is the topological dual of $X$. So, the operator $A$ from (PG)
is (now) a continuous linear operator whose adjoint is
\[
A^{\ast}:Y=\mathbb{R}^{2}\rightarrow X^{\prime}=\mathbb{R}^{\mathbb{N}},\quad
A^{\ast}(u,v)=(u,u+v,...,ku+v,...),\label{r-g3}%
\]
and so (P$_{c,b}^{\ast}$) becomes problem (PG$^{\ast}$).

\medskip At the end of the previous section, we mentioned that Theorems 3.1
and 3.2 from \cite{ViKiTaYe16} provide conditions which ensure that $g(c,b)=0$
when $(c,b)\in\left(  P^{+}+A^{\ast}(Q^{+})\right)  \times\operatorname*{cor}%
\left(  A(P)-Q\right)  $ for the problem (P$_{c,b}$) with $Q=\{0\}$. Besides
such results, in \cite{ViKiTaYe16}, one considers \textquotedblleft the
example of Gale with both $b$ and $c$ being perturbed" to illustrate
\cite[Th.\ 3.1]{ViKiTaYe16}; the problems (P$_{c,b}$) and (P$_{c,b}^{\ast}$)
with $X:=\mathbb{R}^{(\mathbb{N})}$, $Y:=\mathbb{R}^{2}$, $P$ and $Q$ as in
(\ref{r-g2a}), $A$ as in (\ref{r-g2b}) and arbitrary $(c,b)\in X^{\prime
}\times Y$, are denoted by (PG$_{c,b}$) and (PG$_{c,b}^{\ast}$), respectively.

On page 12 of \cite{ViKiTaYe16}, one mentions: \textquotedblleft In connection
with Claim 2, we observe that the question whether the result is valid without
the extra assumption (6) remains open". Our aim in this section is to answer
this open problem.

\medskip In the rest of this section, $F(b)$, $\varphi(c,b)$, $F^{\ast}(c)$,
$\psi(c,b)$, $\Lambda$ and $g$ refer to the problems (PG$_{c,b}$) and
(PG$_{c,b}^{\ast}$). It follows that $A(P)=\{(b_{1},b_{2})\in\mathbb{R}%
^{2}\mid b_{1}\geq b_{2}\geq0\}$ and $c\in P^{+}+\operatorname{Im}A^{\ast}$ if
and only if there exist $(u,v)\in\mathbb{R}^{2}$ and $(\beta_{k}%
)_{k\in\mathbb{N}}\in\mathbb{R}_{+}^{\mathbb{N}}$ such that
\begin{equation}
c_{0}=u+\beta_{0},\ c_{1}=u+v+\beta_{1},\ ...,\
c_{k}=ku+v+\beta_{k},\ \ldots .
\label{r-g4}%
\end{equation}

\begin{lemma}
\label{lem1}Let $c:=(c_{k})_{k\in\mathbb{N}}\in\mathbb{R}^{\mathbb{N}}$,
$(u,v)\in\mathbb{R}^{2}$ and $(\beta_{k})_{k\in\mathbb{N}}\in\mathbb{R}%
_{+}^{\mathbb{N}}$ verify $(\ref{r-g4})$, $\beta^{\prime}\in\lbrack0,\beta
_{0}]$ and $v^{\prime}\in\mathbb{R}$. Then $(u,v)\in F^{\ast}(c)$; moreover
\begin{equation}
(u+\beta^{\prime},v^{\prime})\in F^{\ast}(c)\Leftrightarrow v^{\prime}\leq
v+\inf\{\beta_{k}-\beta^{\prime}k\mid k\geq1\}\Rightarrow\beta^{\prime}%
\leq\overline{\beta}, \label{r-g6}%
\end{equation}
where
\begin{equation}
\overline{\beta}:=\liminf\nolimits_{k\rightarrow\infty}\beta_{k}/k\in
\lbrack0,\infty]. \label{r-g8}%
\end{equation}
Furthermore, if $\beta^{\prime}<\overline{\beta}$, then there exists
$v^{\prime\prime}\in\mathbb{R}$ such that $(u+\beta^{\prime},v^{\prime\prime
})\in F^{\ast}(c)$.
\end{lemma}

Proof. Because (\ref{r-g4}) is verified, $A^{\ast}(u,v)-c\leq0$, and so
$(u,v)\in F^{\ast}(c)$.

Set $v^{\prime\prime}:=v^{\prime}-v$; because $u+\beta^{\prime}\leq c_{0}$,
one has that%
\begin{align*}
(u+\beta^{\prime},v^{\prime}) \in F^{\ast}(c)  &  \Leftrightarrow
\lbrack\forall k\geq1:k(u+\beta^{\prime})+v^{\prime}\leq c_{k}]\Leftrightarrow
\lbrack\forall k\geq1:\beta_{k}\geq\beta^{\prime}k+v^{\prime\prime}]\\
&  \Leftrightarrow v^{\prime}\leq v+\inf\{\beta_{k}-\beta^{\prime}k\mid
k\geq1\},
\end{align*}
whence the equivalence from (\ref{r-g6}) follows. Assume now that $\beta
_{k}\geq\beta^{\prime}k+v^{\prime\prime}$ for $k\geq1$; then $\beta_{k}%
/k\geq\beta^{\prime}+v^{\prime\prime}/k$ for $k\geq1$, whence
\[
\overline{\beta}=\liminf\nolimits_{k\rightarrow\infty}\beta_{k}/k\geq
\liminf\nolimits_{k\rightarrow\infty}(\beta^{\prime}+v^{\prime\prime}%
/k)=\beta^{\prime},
\]
and so the (last) implication from (\ref{r-g6}) holds, too.

Assume that $\beta^{\prime}<\overline{\beta}$; then there exists $k_{0}\geq1$
such that $\beta^{\prime}\leq\beta_{k}/k$ for $k\geq k_{0}$, whence
$k(u+\beta^{\prime})+v=ku+v+k\beta^{\prime}\leq k\overline{u}+v+\beta
_{k}=c_{k}$ for $k\geq k_{0}$. Setting
\[
v_{0}:=\min\{c_{k}-k(u+\beta^{\prime})\mid k\in\overline{1,k_{0}}\},\quad
v^{\prime\prime}:=\min\{v,v_{0}\},
\]
one has $k(u+\beta^{\prime})+v^{\prime\prime}\leq c_{k}$ for $k\geq1$; because
$\beta^{\prime}\leq\beta_{0}$ and $v^{\prime\prime}\leq v$ we get
$k(u+\beta^{\prime})+v^{\prime\prime}\leq c_{k}$ for $k\geq0$, whence
$(u+\beta^{\prime},v^{\prime\prime})\in F^{\ast}(c)$. \hfill$\square$

\medskip

Claim 2 from \cite{ViKiTaYe16} asserts, equivalently, the following:

\textit{Suppose that $b=(b_{1},0)$ with $b_{1}>0$ and $c\in P^{+}%
+\operatorname{Im}A^{\ast}$ is of the form $(\ref{r-g4})$ such that
$\overline{\beta}=0$, where $\overline{\beta}$ is defined in
$(\ref{r-g8})$. Then, $g(c,b)>0$ if $\beta_{0}>0$ and $g(c,b)=0$ if
$\beta_{0}=0$.}

\medskip As recalled above, on page 12 of \cite{ViKiTaYe16} it is said:
\textquotedblleft In connection with Claim 2, we observe that the question
whether the result is valid without the extra assumption (6) remains
open".\footnote{Assumption (6) from \cite{ViKiTaYe16} is equivalent to
$\overline{\beta}=0$, where $\overline{\beta}$ is defined in (\ref{r-g8}).}

\medskip

Related to \cite[Claim 2]{ViKiTaYe16} we have the following result:

\begin{proposition}
\label{p1}Consider $b:=(b_{1},0)$ with $b_{1}>0$ and
\[
c=(u+\beta_{0},u+v+\beta_{1},...,ku+v+\beta_{k},...)\in\mathbb{R}^{\mathbb{N}}
\label{r-g7}%
\]
with $(u,v)\in\mathbb{R}^{2}$ and $(\beta_{k})_{k\in\mathbb{N}}\in
\mathbb{R}_{+}^{\mathbb{N}}$; hence $(c,b)\in\Lambda$. Then $F(b)=\{b_{1}%
e_{0}\}$, and
\begin{equation}
\varphi(c,b)=b_{1}c_{0},\quad\psi(c,b)=b_{1}\overline{u},\quad g(c,b)=b_{1}%
\max\{0,\beta_{0}-\overline{\beta}\}, \label{r-R}%
\end{equation}
where $\overline{\beta}$ is defined in $(\ref{r-g8})$ and $\overline
{u}:=u+\min\{\beta_{0},\overline{\beta}\}$; moreover, $(u,v)$ is an optimal
solution of \emph{(PG}$_{c,b}^{\ast}$\emph{)} if $\min\{\beta_{0}%
,\overline{\beta}\}=0$, and \emph{(PG}$_{c,b}^{\ast}$\emph{)} has optimal
solutions if $0<\beta_{0}<\overline{\beta}$. Furthermore, if $\beta_{0}%
\geq\overline{\beta}>0$, \emph{(PG}$_{c,b}^{\ast}$\emph{)} has optimal
solutions if and only if $\inf\{\beta_{k}-\overline{\beta}k\mid k\geq
1\}\in\mathbb{R}$.
\end{proposition}

Proof. Observe first that $F(b)=\{b_{1}e_{0}\}$, whence $\varphi
(c,b)=b_{1}c_{0}$, and $(u,v)\in F^{\ast}(c)$.

Consider $(u^{\prime},v^{\prime})\in F^{\ast}(c)$; then $u^{\prime}\leq
c_{0}=u+\beta_{0}$ and $ku^{\prime}+v^{\prime}\leq c_{k}=ku+v+\beta_{k}$ for
$k\geq1$. It follows that $u^{\prime}+v^{\prime}/k\leq u+v/k+\beta_{k}/k$ for
$k\geq1$, whence $u^{\prime}\leq u+\overline{\beta}$, where $\overline{\beta}$
is defined in (\ref{r-g8}). It follows that $u^{\prime}\leq\overline{u}\leq
c_{0}$ for all $(u^{\prime},v^{\prime})\in F^{\ast}(c)$.

Because $\psi(c,b)=\sup\{b_{1}u^{\prime}\mid(u^{\prime},v^{\prime})\in
F^{\ast}(c)\}$ and $(u,v)\in F^{\ast}(c)$, we get
\begin{equation}
b_{1}u\leq\psi(c,b)\leq b_{1}\overline{u}\leq b_{1}c_{0}=\varphi(c,b).
\label{r-g10}%
\end{equation}

Assume first that $\min\{\beta_{0},\overline{\beta}\}=0$; then $u=\overline
{u}$, whence $\psi(c,b)=b_{1}u=b_{1}\overline{u}$ by (\ref{r-g10}), and so
$(u,v)$ is an optimal solution of (PG$_{c,b}^{\ast}$).

Assume now that $\min\{\beta_{0},\overline{\beta}\}>0$; take $0<\beta^{\prime
}<\min\{\beta_{0},\overline{\beta}\}$. Using the last part of Lemma
\ref{lem1}, it follows that there exists $v^{\prime\prime}\in\mathbb{R}$ such
that $(u+\beta^{\prime},v^{\prime\prime})\in F^{\ast}(c)$, and so
$\psi(c,b)\geq b_{1}(u+\beta^{\prime})$. Because $\beta^{\prime}$ is arbitrary
in $]0,\min\{\beta_{0},\overline{\beta}\}[$, one obtains that $\psi(c,b)\geq
b_{1}\overline{u}$. Having in view (\ref{r-g10}), one gets $\psi
(c,b)=b_{1}\overline{u}$. Hence $\varphi(c,b)=b_{1}c_{0}$ and $\psi
(c,b)=b_{1}\overline{u}$, and so
\[
g(c,b):=\varphi(c,b)-\psi(c,b)=b_{1}\left(  c_{0}-\overline{u}\right)
=b_{1}\left(  \beta_{0}-\min\{\beta_{0},\overline{\beta}\}\right)  =b_{1}%
\max\{0,\beta_{0}-\overline{\beta}\};
\]
therefore, (\ref{r-R}) holds.

As seen above, $(u,v)$ is an optimal solution of (PG$_{c,b}^{\ast}$) when
$\min\{\beta_{0},\overline{\beta}\}=0$. Assume that $0<\beta_{0}%
<\overline{\beta}$. Using again the last part of Lemma \ref{lem1} with
$\beta_{0}$ instead of $\beta^{\prime}$, one gets $v^{\prime\prime}%
\in\mathbb{R}$ such that $F^{\ast}(c)\ni(u+\beta_{0},v^{\prime\prime
})=(\overline{u},v^{\prime\prime})$, and so $(\overline{u},v^{\prime\prime})$
is an optimal solution of (PG$_{c,b}^{\ast}$).

Assume that $\beta_{0}\geq\overline{\beta}>0$. Because $\psi(c,b)=b_{1}%
\overline{u}=b_{1}(u+\overline{\beta})$, (PG$_{c,b}^{\ast}$) has optimal
solutions if and only if there exists $v^{\prime}\in\mathbb{R}$ such that
$(\overline{u}+\overline{\beta},v^{\prime})\in F^{\ast}(c)$, and this is
equivalent with $\inf\{\beta_{k}-\overline{\beta}k\mid k\geq1\}\in\mathbb{R}$
by the equivalence in (\ref{r-g6}). \hfill$\square$

\medskip

Claim 4 from \cite{ViKiTaYe16} asserts, equivalently, the following:

\textit{For every $c\in P^{+}+\operatorname{Im}A^{\ast}$ and $b=(b_{1},b_{2})$
with $b_{1}>b_{2}>0$, one has $g(c,b)=0$.}

\smallskip The proof of Claim 4 is quite involved in \cite{ViKiTaYe16}; for
its proof one uses \cite[Th.\ 3.1]{ViKiTaYe16}. In fact the conclusion of this
assertion follows immediately using Corollary \ref{c1}. Indeed, because
$A(P)=A(P)-Q=\{(y_{1},y_{2})\in\mathbb{R}^{2}\mid y_{1}\geq y_{2}\geq0\}$,
$\operatorname*{cor}A(P)=\{(y_{1},y_{2})\in\mathbb{R}^{2}\mid y_{1}>y_{2}%
>0\}$, and so condition (iv) of Proposition \ref{p3} is verified; as observed
in Remark \ref{rem1}, $g(c,b)=0$ and Sol(PG$_{c,b}^{\ast}$)$\neq\emptyset$ for
all $(c,b)\in\left(  P^{+}+A^{\ast}(Q^{+})\right)  \times\operatorname*{icr}%
A(P)$.

\medskip

In \cite{KhMoTr19}, only $b$ is perturbed in Gale's example; hence $c:=e_{0}$.
So, one may take $u:=v:=0$, $\beta_{0}:=1$ and $\beta_{k}:=0$ for $k\geq1$ in
(\ref{r-g4}), and so $\overline{\beta}=0$ in (\ref{r-g8}). The results on
problem (PG$_{e_{0},b}$) and (PG$_{e_{0},b}^{\ast}$) from \cite[Claims 2,
3]{KhMoTr19} confirm those obtained in \cite{ViKiTaYe16} for this case; even
more, it is obtained that both problems have optimal solutions for $b\in A(P)$.

\section{On Gale's example from Van Slyke and Wets paper \cite{VanWet68}}

At the end of Section 3.8 of Anderson~\& Nash book \cite{AndNas87} it is
mentioned that \textquotedblleft The example of Section 3.5.1 is due to Gale
and appears in a paper by Van Slyke and Wets (1968)."\footnote{In fact, it is
Section 3.4.1, pages 42, 43.}

Indeed, in Van Slyke and Wets paper \cite{VanWet68}, one finds the following
text containing also the problem to which Gale's example refers:

\smallskip\textquotedblleft Our purpose is to obtain characterizations of the
optimal solutions to

\smallskip Find inf $z=c(x)$ subject to $\left\langle A,x\right\rangle =b$,
$x\in X\subset\mathcal{X}$ (3.1)

\smallskip\noindent where $c(x)$ is a convex functional defined everywhere but
possibly with values $+\infty$ or $-\infty$ for some $x$ in $\mathcal{X}$, $A$
is a continuous linear operator from $\mathcal{X}$ into $\mathcal{Y}$, both
locally convex linear Hausdorf topological spaces and $X$ is a convex subset
of $\mathcal{X}$. By $\left\langle \cdot,\cdot\right\rangle $ we denote linear
composition. It is easy to see that (3.1) is equivalent to

\smallskip Find inf $\eta$ such that $(\eta,y)\in\mathcal{C}\cap
\mathcal{L}\subset\mathbb{R}\times\mathcal{Y}$, \ \ (3.2)

\smallskip\noindent where $\mathbb{R}$ denotes the real line (with its natural topology),

\smallskip$\mathcal{C}=\{(\eta,y)\mid\eta\geq c(x)$, $y=b-\left\langle
A,x\right\rangle $ for some $x\in X\}$ and

\smallskip$\mathcal{L=\{(\eta},0)\mid\eta\in\mathbb{R}$, $0\in\mathcal{Y}\}$.

\smallskip... the following simple example, due to David Gale (which is a
variant of the one given in [9]\footnote{This is our reference \cite{DufKar65}%
.}).

\smallskip

EXAMPLE (3.5). Consider the semi-infinite program

\smallskip Find inf $x_{0}$ subject to $x_{0}+\sum_{n=1}^{\infty}%
nx_{n}=1,\ \sum_{n=1}^{\infty}x_{n}=0,\ x_{n}\geq0$, $n=0,1,..$, and
$\mathcal{X}=\ell_{p}$, $1\leq p<\infty$.

\smallskip

\smallskip\noindent Now

\smallskip$\mathcal{C}=\{(\eta,y_{1},y_{2})\mid\eta\geq x_{0}$, $y_{1}%
=1-x_{0}-\sum_{n=1}^{\infty}nx_{n}$, $y_{2}=0-\sum_{n=1}^{\infty}x_{n}$,
$x_{n}\geq0$, $n=0,1,...\}$

\smallskip We first observe that

\smallskip$\overline{\mathcal{L}\cap\mathcal{C}}=\mathcal{L}\cap
\mathcal{C}=\{(\eta,0,0)\mid\eta\geq1\}$.

\smallskip\noindent But since the closure of $\mathcal{C}$ is

\smallskip$\overline{\mathcal{C}}=\{(\eta,y_{1},y_{2})\mid\eta\geq x_{0}$,
$y_{1}=1-x_{0}-\sum_{n=1}^{\infty}nx_{n}-z$, $y_{2}=0-\sum_{n=1}^{\infty}%
x_{n}$, $x_{n}\geq0$, $n=0,1,...$, $z\geq0\}$,

\smallskip\noindent we have that $\mathcal{L}\cap\overline{\mathcal{C}%
}=\{(\eta,0,0)\mid\eta\geq0\}$. Note that in this example $\mathcal{L}%
\cap\mathcal{C}$ is closed and thus the inf is attained."

\medskip

From the above text, we understand that $\mathcal{X}$, $\mathcal{Y}$ are real
H.l.c.s., $A:\mathcal{X}\rightarrow\mathcal{Y}$ is a continuous linear
operator, $c:\mathcal{X}\rightarrow\overline{\mathbb{R}}$ is convex and
$X\subset\mathcal{X}$ is a convex set; we set $Ax$ instead of $\left\langle
A,x\right\rangle $. Moreover,%
\[
\mathcal{C}=\{(\eta,b-Ax\mid x\in X,\ \eta\in\mathbb{R}\cap\lbrack
c(x),\infty\lbrack\}\subset\mathbb{R\times}\mathcal{Y},\quad\mathcal{L}%
=\mathbb{R}\times\{0\}\subset\mathbb{R\times}\mathcal{Y}. \label{r-vw1}%
\]

Below, we refer to Example (3.5) from \cite{VanWet68} without mentioning it
explicitly. Having in view the formulation of this example, in the sequel, all
the sequences are indexed by $n\in\mathbb{N}$; for example, by $x:=(x_{n})$ we
mean $x:=(x_{n})_{n\in\mathbb{N}}$ $(\in\mathbb{R}^{\mathbb{N}})$.

So, $\mathcal{X}:=\ell_{p}$ with $1\leq p<\infty$ and $\mathcal{Y}%
:=\mathbb{R}^{2}$; consequently, $\mathcal{X}^{\ast}$ is identified with
$\ell_{q}$, where $q:=p/(p-1)\in{}]1,\infty]$ for $p\in{}]1,\infty\lbrack$,
$q:=\infty$ for $p=1$, and $\left\langle x,y\right\rangle =\sum_{n\geq0}%
x_{n}y_{n}$ for $x:=(x_{n})\in\mathcal{X}$ and $y:=(y_{n})\in\mathcal{X}%
^{\ast}$. Moreover, $\mathcal{Y}^{\ast}$ is identified with $\mathbb{R}^{2}$
$(=\mathcal{Y})$.

\smallskip For $z:=(z_{n})\in\mathbb{R}^{\mathbb{N}}$, by ${\textstyle\sum
\nolimits_{n=0}^{\infty}}z_{n}\in\mathbb{R}$ we mean that the series
${\textstyle\sum\nolimits_{n=0}^{\infty}}z_{n}$ is convergent in $\mathbb{R}$
(and so its sum belongs to $\mathbb{R}$ and $z\in\mathfrak{c}_{0}%
:=\{x\in\mathbb{R}^{\mathbb{N}}\mid x_{n}\rightarrow0\}$).

\begin{remark}
\label{rem-vw2}Consider the sets
\begin{align*}
D_{0}^{a}  &  :=\{z\in\mathbb{R}^{\mathbb{N}}\mid{\textstyle\sum
\nolimits_{n=0}^{\infty}}n\left\vert z_{n}\right\vert \in\mathbb{R}%
\},\ \ D_{0}:=\{z\in\mathbb{R}^{\mathbb{N}}\mid{\textstyle\sum\nolimits_{n=0}%
^{\infty}}nz_{n}\in\mathbb{R}\},\\
D_{1}  &  :=\{z\in\mathbb{R}^{\mathbb{N}}\mid{\textstyle\sum\nolimits_{n=0}%
^{\infty}}z_{n}\in\mathbb{R}\}.
\end{align*}
Then $D_{0}^{a}$, $D_{0}$ and $D_{1}$ are linear spaces such that
$\mathbb{R}^{(\mathbb{N})}\subset D_{0}^{a}\subset D_{0}\cap\ell_{1}%
\subset\mathfrak{c}_{0}$ and $D_{0}\subset D_{1}\cap\ell_{r}$ for every
$r\in{}]1,\infty\lbrack$, each inclusion being strict; moreover,
$D_{0}\not \subset \ell_{1}$.
\end{remark}

Proof. To check that $D_{0}^{a}$, $D_{0}$ and $D_{1}$ are linear spaces is routine.

The inclusions $\mathbb{R}^{(\mathbb{N})}\subset D_{0}^{a}\subset D_{0}%
\cap\ell_{1}\subset\mathfrak{c}_{0}$ are obvious; moreover, $(2^{-n})\in
D_{0}^{a}\setminus\mathbb{R}^{(\mathbb{N})}$, $\left(  \frac{(-1)^{n}%
}{(n+1)^{2}}\right)  \in(D_{0}\cap\ell_{1})\setminus D_{0}^{a}$ and $\left(
\frac{1}{n+1}\right)  \in\mathfrak{c}_{0}\setminus\ell_{1}$.

The inclusion $D_{0}\subset D_{1}$ follows by Abel's test. Take $z:=(z_{n})\in
D_{0}$ and $r\in{}]1,\infty\lbrack$. Then $(nz_{n})$ is bounded, whence there
exists $\mu>0$ such that $\left\vert nz_{n}\right\vert \leq\mu$ for
$n\in\mathbb{N}$, and so $\left\vert z_{n}\right\vert ^{r}\leq\mu^{r}/n^{r}$
for $n\geq1$; hence ${\textstyle\sum\nolimits_{n=0}^{\infty}} \left\vert
z_{n}\right\vert ^{r}$ is convergent. Moreover, $\left(  \frac{(-1)^{n}}%
{n+1}\right)  \in(D_{1}\cap\ell_{r})\setminus D_{0}$ for $r\in{}%
]1,\infty\lbrack$. Furthermore, take $z:=(z_{n})$ defined by $z_{0}:=z_{1}%
:=0$, $z_{n}:=\frac{(-1)^{n}}{n\ln n}$ for $n\geq2$; then $\sum_{n\geq0}%
nz_{n}=\sum_{n\geq2}\frac{(-1)^{n}}{\ln n}\in\mathbb{R}$ by Leibniz's test,
but $\sum_{n\geq0}\left\vert z_{n}\right\vert =\sum_{n\geq2}\frac{1}{n\ln
n}=\infty$, and so $z\in D_{0}\setminus\ell_{1}$. \hfill$\square$

\medskip Having in view the description of the data considered in (3.1), in
Example\ (3.5), for $x:=(x_{n})$, one must have that $c(x):=x_{0}$,
\begin{equation}
Ax:=\left(  x_{0}+{\sum\nolimits_{n=1}^{\infty}}nx_{n},{\sum\nolimits_{n=1}%
^{\infty}}x_{n}\right)  \in\mathbb{R}^{2}, \label{r-vw2}%
\end{equation}
$b:=(1,0)\in\mathbb{R}^{2}$ and $X:=P:=(\ell_{p})_{+}:=\ell_{p}\cap
\mathbb{R}_{+}^{\mathbb{N}}$. It follows that $c(x)=\left\langle
x,e_{0}\right\rangle $ for $x\in\mathcal{X}$, and $P^{+}=(\ell_{q})_{+}$.

\begin{remark}
\label{rem-vw3}Let $D_{A}:=\{x\in\mathcal{X}\mid Ax\in\mathcal{Y}\}$ be the
(maximal) domain of $A$ as an operator from $\mathcal{X}$ to $\mathcal{Y}$.
Then $D_{A}=\mathcal{X}\cap D_{0}$, $D_{A}$ is a dense linear subspace of
$\mathcal{X}$, and $A:D_{A}\rightarrow\mathcal{Y}$ is a linear operator;
moreover, $\operatorname*{gph}A$ is not closed. The same conclusions hold for
$D_{A}^{a}:=\mathcal{X}\cap D_{0}^{a}$ and $A|_{D_{A}^{a}}$.
\end{remark}

Proof. From Remark \ref{rem-vw2} and (\ref{r-vw2}) we obtain that
$D_{A}=\mathcal{X}\cap D_{0}\cap D_{1}=\mathcal{X}\cap D_{0}\supseteq
\mathbb{R}^{(\mathbb{N})}$, and so $D_{A}$ is a dense linear subspace of
$\mathcal{X}$. From the expression of $Ax$, it follows immediately that $A$ is linear.

Assume that $G:=\operatorname*{gph}A$ $(=\{(x,Ax)\mid x\in D_{A}\})$ is a
closed subset of $\mathcal{X}\times\mathcal{Y}$, and so $G$ is a closed linear
subspace of the normed vector space $\mathcal{X}\times\mathcal{Y}$. Because
$D_{A}\neq\mathcal{X}$, there exists $\overline{x}\in\mathcal{X}\setminus
D_{A}$, and so $\overline{x}\neq0$, $E:=\{\overline{x}\}\times\mathcal{Y}$ is
a closed convex subset of $\mathcal{X}\times\mathcal{Y}$ and $E\cap
G=\emptyset$; moreover, $E$ is (clearly) locally compact because
$\dim(\operatorname*{span}E)=3$. Furthermore, the asymptotic cones of $E$ and
$G$ are $C_{E}:=\{0\}\times\mathcal{Y}$ and $C_{G}:=G$, and so $C_{E}\cap
C_{G}=\{(0,0)\}$. Applying \cite[Prop.\ 1]{Die66}, the set $E-G$
$(\subset\mathcal{X}\times\mathcal{Y})$ is a closed convex set with
$(0,0)\notin E-G$. Therefore, there exist $(u,v)\in\mathcal{X}^{\ast}%
\times\mathcal{Y}^{\ast}$ and $\gamma\in\mathbb{R}$ such that
\[
0>\gamma\geq\left\langle (\overline{x}-x,u\right\rangle +\left\langle
y-Ax,v\right\rangle =\left\langle \overline{x},u\right\rangle +\left\langle
y,v\right\rangle -\left\langle x,u\right\rangle -\left\langle
Ax,v\right\rangle \quad\forall x\in D_{A},\ y\in\mathcal{Y}.
\]
Because $D_{A}$ and $\mathcal{Y}$ are linear spaces, one obtains that
$\left\langle y,v\right\rangle =0$ for $y\in\mathcal{Y}$ and $\left\langle
x,u\right\rangle +\left\langle Ax,v\right\rangle =0$ for $x\in D_{A}$, and so
$\left\langle \overline{x},u\right\rangle <0$, $v=0$ and $\left\langle
x,u\right\rangle =0$ for every $x\in D_{A}$. Since $D_{A}$ is dense in
$\mathcal{X}$ and $u\in\mathcal{X}^{\ast}$, we get $u=0$, contradicting the
fact that $\left\langle \overline{x},u\right\rangle <0$. Hence
$\operatorname*{gph}A$ is not closed. \hfill$\square$

\medskip Remark \ref{rem-vw3} shows that Example (3.5) is not adequate for
illustrating the problem (3.1).

\smallskip Having in view Remark \ref{rem-vw3}, the definition of
$\mathcal{C}$ in Example\ (3.5) is not clear (because $\sum_{n=1}^{\infty
}x_{n}$ and $\sum_{n=1}^{\infty}nx_{n}$ do not make sense for every
$x\in\mathcal{X}$). Let us consider that $\mathcal{C}$ is the set%
\[
\mathcal{C}^{\prime}=\{(\eta,y)\in\mathbb{R}\times\mathcal{Y}\mid\exists
x\in\mathbb{R}_{+}^{\mathbb{N}}:\eta\in\lbrack x_{0},\infty\lbrack
,\ y_{1}=1-x_{0}- {\textstyle\sum\nolimits_{n=1}^{\infty}} nx_{n}\in
\mathbb{R},\ y_{2}=- {\textstyle\sum\nolimits_{n=1}^{\infty}} x_{n}%
\in\mathbb{R}\},
\]
that is, $\mathcal{C}^{\prime}$ is the biggest set for which all the elements
from the description of $\mathcal{C}$ are well defined;\footnote{This agrees
with the following text from \cite{DufKar65}: \textquotedblleft%
\emph{Definition 2}. The infinite program P is said to be \emph{consistent} if
there exists a sequence $(x_{1},x_{2},...)$ for which all of the series in
question converge and which satisfies all of the constraints. Such a sequence
is termed \emph{feasible}.".} hence $x$ from the definition of $\mathcal{C}%
^{\prime}$ belongs to $\mathbb{R}_{+}^{\mathbb{N}}\cap D_{A}$.

Let $(\eta,y)\in\mathcal{C}^{\prime}$, and take $x\in\mathbb{R}_{+}%
^{\mathbb{N}}$ such that $\eta\in\lbrack x_{0},\infty\lbrack$, $y_{1}%
=1-x_{0}-\sum_{n=1}^{\infty}nx_{n}\in\mathbb{R}$ and$\ y_{2}=-\sum
_{n=1}^{\infty}x_{n}\in\mathbb{R}$. Then $x\in\ell_{1}$ because $x_{n}\geq0$
for $n\geq0$ and $\sum_{n=1}^{\infty}x_{n}\in\mathbb{R}$, and so $x\in\ell
_{r}$ for $r\in\lbrack1,\infty\lbrack$; in particular, $x\in\mathcal{X}$. It
follows that $y_{2}\leq0$ and
\[
y_{1}=1-x_{0}+y_{2}-{\textstyle\sum\nolimits_{n=2}^{\infty}} (n-1)x_{n}%
\leq1-x_{0}+y_{2}\leq1+y_{2},
\]
whence $y_{1}\leq1+y_{2}\leq1$ and $x_{0}\leq1+y_{2}-y_{1}$; moreover, if
$y_{2}=0$, then $x_{n}=0$ for $n\geq1$ and so $y_{1}=1-x_{0}$, whence
$\eta\geq x_{0}=1-y_{1}$. Therefore, $\mathcal{C}^{\prime}\subset
\mathcal{C}_{0}\cup\mathcal{C}_{1}$, where
\begin{equation}
\mathcal{C}_{0}:=\{(\eta,y_{1},0)\mid y_{1}\leq1,\ \eta\geq1-y_{1}%
\},\quad\mathcal{C}_{1}:=\{(\eta,y_{1},y_{2})\mid y_{1}-1\leq y_{2}%
<0,\ \eta\geq0\}. \label{r-vw3}%
\end{equation}

Conversely, take first $(\eta,y_{1},y_{2})\in\mathcal{C}_{0}$. Then $y_{2}=0$,
$y_{1}\leq1$ and $\eta\geq1-y_{1}$ $(\geq0)$, and so $x:=(1-y_{1},0,0,...)$
verifies the conditions from the definition of $\mathcal{C}^{\prime}$; hence
$(\eta,y_{1},y_{2})\in\mathcal{C}^{\prime}$. Take now $(\eta,y_{1},y_{2}%
)\in\mathcal{C}_{1}$; then $\alpha\geq\beta>0$, where $\alpha:=1-y_{1}$ and
$\beta:=-y_{2}>0$. Take $\overline{n}\in\mathbb{N}$ with $\overline{n}\geq
\max\{2,\alpha/\beta\}$, $x_{1}:=(\overline{n}\beta-\alpha)/(\overline{n}-1)$,
$x_{\overline{n}}:=(\alpha-\beta)/(\overline{n}-1)$ and $x_{n}:=0$ for
$n\in\mathbb{N}\setminus\{1,\overline{n}\}$; then this $x$ verifies the
conditions from the definition of $\mathcal{C}^{\prime}$, and so $(\eta
,y_{1},y_{2})\in\mathcal{C}^{\prime}$. Therefore, $\mathcal{C}^{\prime
}=\mathcal{C}_{0}\cup\mathcal{C}_{1}$, where $\mathcal{C}_{0}$ and
$\mathcal{C}_{1}$ are defined in (\ref{r-vw3}).

\medskip

Having in view the description of $\mathcal{C}$ and $\overline{\mathcal{C}}$
in \cite[Example\ (3.5)]{VanWet68}, in the case in which these sets were
correctly defined, one would have $\overline{\mathcal{C}}=\mathcal{C}%
-(\{0\}\times\mathbb{R}_{+}\times\{0\})$. However, we have that
\[
\operatorname*{cl}\mathcal{C}^{\prime}=\{(\eta,y_{1},y_{2})\mid y_{1}-1\leq
y_{2}\leq0,\ \eta\geq0\}.
\]
Consequently, $\overline{\mathcal{L}\cap\mathcal{C}^{\prime}}=\mathcal{L}%
\cap\mathcal{C}^{\prime}=\{(\eta,0,0)\mid\eta\geq1\}$ and $\mathcal{L}%
\cap\overline{\mathcal{C}^{\prime}}=\{(\eta,0,0)\mid\eta\geq0\}$.

\medskip

In conclusion, Example\ (3.5) is not adequate for illustrating the problem
(3.1) from \cite{VanWet68} because the operator $A$ is not a continuous linear
operator from $\mathcal{X}$ into $\mathcal{Y}$; moreover, because the graph of
$A$ is not closed, one can not speak of its adjoint, which is necessary when
defining (using) the dual problem. However, the sets $\overline{\mathcal{L}%
\cap\mathcal{C}}$ and $\mathcal{L}\cap\overline{\mathcal{C}}$ are those
indicated in \cite[ Example\ (3.5)]{VanWet68}.

\end{document}